NOTE ON A COIN TOSSING PROBLEM POSED BY DANIEL LITT

**Bruce Levin** BL6@COLUMBIA.EDU
*Department of Biostatistics*
*Mailman School of Public Health*
*Columbia University*
*New York, NY 10032, USA*

**Abstract**

We present an analysis of a coin-tossing problem posed by Daniel Litt which has generated some popular interest. We demonstrate a recursive identity which leads to relatively simple formulas for the excess number of wins for one player over the other together with its increments as the number of coin tosses increases. Formulas and recursive algorithms are provided to calculate the number of sequences with any given point-score difference.

Daniel Litt of the Uniiversity of Toronto posed the following coin tossing problem, discussed in E. Klarreich's article, Perplexing the Web, One Probability Puzzle at a Time (Klarreich, 2024):

> Alice and Bob flip a (fair) coin 100 times. Anytime there are two heads in a row, Alice gets a point; when a head is followed by a tail, Bob gets a point. So in the sequence THHHT, Alice gets two points and Bob gets one. Who is more likely to win?

The purpose of this note is to prove that for any version of the game with $n>2$ tosses, Bob has more winning binary $n$-sequences than Alice, hence Bob is more likely to win for any such $n$. The approach relies on some interesting recursive identities among "heady close-call" binary sequences, defined below.

**1. Definitions and notation.**

Let $x^{(n)}=(x_1^{(n)},...,x_n^{(n)})$ be a binary sequence of length $n$ with 1's representing heads and 0's tails. Let $S(x^{(n)})$ denote Alice's point total minus Bob's point total given sequence $x^{(n)}$; in symbols,

$$S(x^{(n)})=\sum_{i=1}^{n-1} I[x_i^{(n)}=x_{i+1}^{(n)}=1]-\sum_{i=1}^{n-1} I[x_i^{(n)}=1,\, x_{i+1}^{(n)}=0].$$

A win for Alice (henceforth $A$) occurs when $S(x^{(n)})>0$ and a win for Bob (henceforth B) occurs when $S(x^{(n)})<0$. Ties with $S(x^{(n)})=0$ can occur but are not counted as wins for either $A$ or $B$.

Let $D_n$ be the total number of winning $n$-sequences for *B* *minus* the number of winning $n$-sequences for $A$,

$$D_n = \sum_{x^{(n)}} I[S(x^{(n)}) < 0] - \sum_{x^{(n)}} I[S(x^{(n)}) > 0],$$

where the sums are over all $2^n$ binary sequences $x^{(n)}$. Also define the *forward increment from n to* $n+1$ to be $\Delta_{n+1} = D_{n+1} - D_n$ for $n \geq 2$. We wish to show $D_n > 0$ for all $n > 2$ and our approach will be to show that $\Delta_{n+1} > 0$ for all $n \geq 2$.

**2. Proof of positive increments.**

To analyze the game it will be helpful to imagine cross-classifying all $2^n$ binary sequences by the outcome of the final toss, $x_n^{(n)} = 1$ vs. $x_n^{(n)} = 0$, versus a relevant five-category classification of $S(x^{(n)})$, namely, $S(x^{(n)}) > 1$, $S(x^{(n)}) = 1$, $S(x^{(n)}) = 0$, $S(x^{(n)}) = -1$, and $S(x^{(n)}) < -1$. We denote the frequencies in the first row of the resulting 2×5 table by $h_1(n),\ldots,h_5(n)$ and those in the second row by $t_1(n),\ldots,t_5(n)$, with the mnemonic *h* or *t* referring to $x_n^{(n)} = 1$ and $x_n^{(n)} = 0$, respectively. The notation is summarized in the diagram below. When there is no risk of ambiguity we may omit ($n$) from the notation.

| Last toss | Total points for *A* minus total points for *B* | | | | | Total |
|---|---|---|---|---|---|---|
| | $S(x) > 1$ | $S(x) = 1$ | $S(x) = 0$ | $S(x) = -1$ | $S(x) < -1$ | |
| $x_n = 1$ | $h_1$ | $h_2$ | $h_3$ | $h_4$ | $h_5$ | $2^{n-1}$ |
| $x_n = 0$ | $t_1$ | $t_2$ | $t_3$ | $t_4$ | $t_5$ | $2^{n-1}$ |
| Total | $h_1 + t_1$ | $h_2 + t_2$ | $h_3 + t_3$ | $h_4 + t_4$ | $h_5 + t_5$ | $2^n$ |

Then

$$D_n = \{h_4(n) + t_4(n) + h_5(n) + t_5(n)\} - \{h_1(n) + t_1(n) + h_2(n) + t_2(n)\}.$$

Let us call sequences $x^{(n)}$ with $S(x^{(n)}) = \pm 1$ and $x_n^{(n)} = 1$ *heady close-call winning sequences* or *heady close-call wins*. The first remarkable fact about the $D_n$ sequence is that the the increments $\Delta_{n+1}$ depend only on $D_n$ and (at most) the numbers of heady close-call wins for *A* and for *B*. This is Lemma 1.

Lemma 1. $\Delta_{n+1} = D_n - \{h_4(n) - h_2(n)\}$ for all $n \geq 2$.

Proof of Lemma 1. We track the contribution that each $x^{(n)}$ makes to the change in the win-count difference to $D_{n+1}$ from $D_n$ as an additional toss is made. We need only take note of $x^{(n)}$ that produce a net change in the win-count difference of plus or minus 1 as we add the two possible outcomes of $x_{n+1} = 1$ or $x_{n+1} = 0$, because sequences which don't alter the win counts don't contribute to $\Delta_{n+1}$. For example, a sequence $x^{(n)}$ of type $x_n = 1$ and $S(x^{(n)}) > 1$, which yields a win for *A* after *n* tosses, yields two wins for *A* after $n+1$ tosses, when either $x_{n+1} = 1$ or $x_{n+1} = 0$. There is thus a net increase of –1 in the total number of wins for *B* over *A*. This is recorded as –1 under the column labelled "Net contribution to $\Delta_{n+1}$" in the first line of the chart below. By contrast, a sequence $x^{(n)}$ of type $x_n = 1$ and $S(x^{(n)}) = 1$, which also yields a win for *A* after *n* tosses, continues to yield one win for *A* if $x_{n+1} = 1$ but becomes a tied sequence if $x_{n+1} = 0$. There is thus no net change in the total number of wins for *B* over *A* due to such $x^{(n)}$, i.e., their net contribution to $\Delta_{n+1}$ is 0. Continuing in this way we see from the chart that the only types of sequences from the first row of the 2×5 table which make non-zero contributions to $\Delta_{n+1}$ are from columns 1 and 5, while all but the tied sequence types from the second row do contribute.

| Sequence type $x_n^{(n)}, S(x^{(n)})$ | Frequency | $S(x^{(n)},1)$ | $S(x^{(n)},0)$ | Net contribution to $\Delta_{n+1}$ |
|---|---|---|---|---|
| 1, $>1$ | $h_1$ | $>1$ | $\geq 1$ | –1 |
| 1, $= 1$ | $h_2$ | $>1$ | $= 0$ | 0 |
| 1, $= 0$ | $h_3$ | $= 1$ | $= -1$ | 0 |
| 1, $= -1$ | $h_4$ | $= 0$ | $< -1$ | 0 |
| 1, $<-1$ | $h_5$ | $\leq -1$ | $< -1$ | +1 |
| 0, $>1$ | $t_1$ | $>1$ | $> 1$ | –1 |
| 0, $= 1$ | $t_2$ | $= 1$ | $= 1$ | –1 |
| 0, $= 0$ | $t_3$ | $= 0$ | $= 0$ | 0 |
| 0, $= -1$ | $t_4$ | $= -1$ | $= -1$ | +1 |
| 0, $<-1$ | $t_5$ | $<-1$ | $< -1$ | +1 |

Therefore, summing the contributions to $\Delta_{n+1}$ over all sequence types with their respective frequencies gives

$$\begin{aligned}\Delta_{n+1} &= -h_1(n) + h_5(n) - t_1(n) - t_2(n) + t_4(n) + t_5(n) \\ &= \{h_4(n) + t_4(n) + h_5(n) + t_5(n)\} - \{h_1(n) + t_1(n) + h_2(n) + t_2(n)\} - \{h_4(n) - h_2(n)\} \\ &= D_n - \{h_4(n) - h_2(n)\}\end{aligned}$$

as was to be shown. □

The recursive identity in the next lemma, which we will prove in the next section, will, remarkably, identify $\Delta_{n+1}$ as $h_2(n)$, the number of heady close-call wins for *A*, and $D_n$ as $h_4(n)$, the number of heady close-call wins for *B*.

Lemma 2. $h_4(n+1) = h_2(n) + h_4(n)$ for all $n \geq 2$.

Lemmas 1 and 2 imply the main result, as follows.

Theorem 1. For all $n \geq 2$, (i) $D_n = h_4(n)$ and (ii) $\Delta_{n+1} = h_2(n)$.

Proof. By induction on $n$. For $n=2$, we have the 2×5 table

| 0 | 1 | 1 | 0 | 0 |
|---|---|---|---|---|
| 0 | 0 | 1 | 1 | 0 |
| 0 | 1 | 2 | 1 | 0 |

with $D_2 = 1-1 = 0 = h_4(2)$, which is (i). For $n=3$ we have the 2×5 table

| 1 | 1 | 1 | 1 | 0 |
|---|---|---|---|---|
| 0 | 0 | 2 | 2 | 0 |
| 1 | 1 | 3 | 3 | 0 |

with $D_3 = (3+0)-(1+1) = 1$, so that $\Delta_3 = D_3 - D_2 = 1-0 = 1 = h_2(2)$, which is (ii). So assuming that (i) and (ii) hold up to some $n$, we are to show that they hold for $n+1$. For (i), $D_{n+1} = \Delta_{n+1} + D_n = h_2(n) + h_4(n) = h_4(n+1)$, which is (i) for $n+1$. The first equality is by definition, the second holds by the inductive hypotheses, and the third holds by Lemma 2. For (ii), $\Delta_{n+2} = D_{n+1} - \{h_4(n+1) - h_2(n+1)\} = h_4(n+1) - \{h_4(n+1) - h_2(n+1)\} = h_2(n+1)$, which is (ii) for $n+1$. The first equality holds by Lemma 1 and the second holds by (i) for the case $n+1$, as was just shown. □

It follows that $D_n > 0$ for all $n \geq 3$ because for such $n$, there is always a heady close-call win for $B$, namely $x^{(n)} = (1,0,...,0,1)$, so by Theorem 1, $h_4(n) \geq 1$ whence $D_n = h_4(n) \geq 1$. The exception for $n=2$ where $D_2 = h_4(2) = 0$ arises because there is no available "room" for an interior 0. Furthermore, $\Delta_{n+1} = h_2(n) > 0$ for all $n \geq 2$, because there is always a heady close-call win for $A$, namely, $x^{(2)} = (1,1)$ or $x^{(n)} = (0,...,0,1,1)$ for $n>2$. Thus Bob always has more winning sequences than Alice starting with three tosses and the gap between the number of Bob's and Alice's winning sequences forever widens.

Note that we did not need to evaluate $h_2(n)$ or $h_4(n)$ explicitly for $n>3$ to draw the above conclusions. The proof of Lemma 2 does provide a lovely, explicit formula for $h_2(n)$, so there is no mystery about the growth of $h_4(n)$. We turn to that next.

**3. Proof of Lemma 2 and formulas for heady close-call wins.**

We demonstrate that $h_4(n+1) = h_2(n) + h_4(n)$ for all $n \geq 2$, the proof of which will provide simple formulas for $h_2(n)$ and $h_4(n)$.

Clearly, any sequence $x^{(n)}$ with $x_n^{(n)}=1$ and $S(x^{(n)})=-1$, of which there are $h_4(n)$, generates a sequence $x^{(n+1)}=(0,x^{(n)})$ with $x_{n+1}^{(n+1)}=1$ and $S(x^{(n+1)})=-1$, because leading zeros do not alter the value of $S(x)$. So it will suffice to show there are precisely $h_2(n)$ additional sequences $x^{(n+1)}$ with $x_{n+1}^{(n+1)}=1$, $S(x^{(n+1)})=-1$, and $x_1^{(n+1)}=1$ (else the sequence would already have been counted among those among the first $h_4(n)$).

Definition: Given any sequence $x^{(n)}$, suppose we mark down a + sign each time two consecutive 1's occur or a – sign if a 1,0 occurs in sequence, ignoring the 0,1 or 0,0 pairs. We call the pattern of + and – signs the *signature* of $x^{(n)}$ and denote it by $\sigma=\sigma(x^{(n)})$.

For example, the 8-sequence (0,0,1,1,1,0,0,1), a heady close-call win for *A*, has signature $\sigma=++-$, while (0,1,0,1,0,0,1,1), a heady close-call win for *B*, has signature $\sigma=--+$. This is an example of a complementary signature.

Definition: Given a signature $\sigma$, the *complementary* signature $\overline{\sigma}$ interchanges the + and – signs.

When considering only heady close-call winning sequences, the number of + signs differs from the number of – signs by plus or minus 1. So the signatures arising from heady close call wins for *B* are in one-to-one correspondence with those for *A*, namely, as their complements.

In general, several *n*-sequences can have the same signature. We will show that for any given signature with one more + sign than – sign, the total number of heady close-call winning *n*-sequences for *A* with the given signature *exactly equals* the number of heady close-call winning $(n+1)$-sequences for *B* that begin with a 1 and have the *complementary* signature. Summing over all signatures from heady close-call winning *n*-sequences for *A* provides the required number $h_2(n)$ of $(n+1)$-sequences $x^{(n+1)}$ beginning with a 1 that are heady close-call wins for *B* with the complementary signature. Conversely, any such $(n+1)$-sequence will have a signature that must be the complement of some signature among those from heady close-call winning *n*-sequences for *A*, so $h_2(n)$ is precisely the number of additional $(n+1)$-sequences comprising $h_4(n+1)$.

We establish the desired identity by exhibiting an algorithm that generates all heady close-call sequences of either type having a given signature. The algorithm will generate the same number of sequences in either case.

Definition. For a given signature $\sigma$, the *heady minimum-length* sequence $\mu=\mu(\sigma)$ with that signature specifies a 1,0 pair for each – and a string of consecutive 1's for each string of consecutive +'s (the former one unit longer than the latter). For a string of + signs followed by a – sign, $\mu(\sigma)$ simply appends a 0 after the string of 1's. A final 1 is appended if $\sigma$ ends with a –. A final 1 is already present if $\sigma$ ends with a +. Also, let the length of the heady minimum-length sequence be denoted by $\lambda = \lambda(\mu) = \lambda(\mu(\sigma))$. We may omit the adjective "heady" below but we always intend the last element of $\mu(\sigma)$ to be 1.

For example, given signature $\sigma = ++-$, $\mu(\sigma) = (1,1,1,0,1)$. Given signature $--+$, $\mu(\sigma) =$ (1,0,1,0,1,1).

For a given signature $\sigma$ of a close-call win, let $k = k(\sigma)$ denote the number of + signs in $\sigma$. Then there are *k initial* 1's in the minimum-length sequence, where we count only the *first* 1 in a string of contiguous 1's as an initial 1. Now let $m = m(n,\sigma) = n - \lambda = n - \lambda(\mu(\sigma))$, which gives the total number of 0's that can be inserted immediately in front of initial 1's to comprise a sequence of length *n*. Then a multinomial partition of *m* units into *k* bins will specify how many additional zeros to insert in front of each initial 1. The total number of such partitions equals $\binom{m+k-1}{k-1}$ by a stars-and-bars argument.

For example, given signature $\sigma = ++-$, the minimum-length sequence $\mu(\sigma) = (1,1,1,0,1)$ is of length $\lambda(\sigma) = 5$ with $k$=2 initial 1's. To generate all heady close-call winning sequences for *A* of length $n$=8, say, with the given signature, since $m$=8–5=3, we have $\binom{3+2-1}{2-1} = 4$ partitions of 3 zeros into 2 bins, namely, (3,0), (2,1), (1,2), and (0,3). The first component specifies how many 0's to insert before the first initial 1 and the second specifies how many 0's to insert before the second initial 1. For the partition (2,1), for example, the algorithm outputs the 8-sequence (0,0,1,1,1,0,0,1), while for the partition (0,3), the algorithm outputs (1,1,1,0,0,0,0,1). Thus there are 4 heady close-call 8-sequence wins for *A* with signature $++-$.

As another example, consider the signature $\sigma = +--+-+$ and suppose we wish to generate all heady close-call wins for *A* of length 13 with that signature. The minimum-length sequence is $\mu(\sigma)$ $= (1,1,0,1,1,0,1,1)$ of length $\lambda(\mu) = 8$ with $k$=3 initial 1's, allowing $m$=13–8=5 zeros to insert. Then there are $\binom{6+3-1}{3-1} = \binom{7}{2} = 21$ trinomial partitions of 5 into 3 bins. For the partition (2, 1, 2), for example, the algorithm outputs the 13-sequence (0,0,1,1,0,0,1,1,0,0,0,1,1).

Now consider generating all heady close-call winning $(n+1)$-sequences for *B* starting and ending with 1 with the complementary signature $\overline{\sigma}$. We again obtain the minimum-length sequence $\mu(\overline{\sigma})$, adding a 1 at the end if $\overline{\sigma}$ ends with –. For this sequence type, the minimal-length sequence will always be one unit longer than that of the original signature, $\lambda(\mu(\overline{\sigma})) = 1 + \lambda(\mu(\sigma))$, so that $m(n+1,\overline{\sigma}) = n+1-\lambda(\mu(\overline{\sigma})) = n+1-1-\lambda(\mu(\sigma)) = n-\lambda(\mu(\sigma)) = m(n,\sigma)$, i.e., we have the same number of excess 0's to insert as for *n*-sequences with the original signature. Now, however, we *do not allow* any 0's to be inserted in front of the automatic leading 1 of the $(n+1)$-sequence, so that $k = k(\sigma)$ of the original signature still counts the number of initial 1's in front of which to insert 0's. *Therefore the algorithm generates exactly the same number of multinomial partitions by inserting the corresponding number of zeros in front of the other initial-1 positions.*

In the above example, the complementary signature is $\overline{\sigma} = -+-+-$ with minimal-length sequence $\mu(\overline{\sigma}) =$ (1,0,1,1,0,1,1,0,1) of length $\lambda(\mu(\overline{\sigma})) = 9$ with $k$=3 initial-1 positions (ignoring the automatic leading 1) and $m$=14–9=5 as before. The same 21 trinomial partitions of 5 zeros into 3 bins generate all the heady close-call winning sequences for *B* with the given complementary signature and starting with a 1. For example, the partition (2, 1, 2) now outputs the 14-sequence (1,0,0,0,1,1,0,0,1,1,0,0,0,1). This concludes the proof of Lemma 2. □

From the above one-to-one correspondences, we get the following useful formula for $h_2(n)$.

<u>Corollary to Lemma 2</u>: $h_2(n) = \sum_{k=1}^{\lfloor (n+1)/3 \rfloor} \binom{2k-1}{k}\binom{n-2k}{k-1}$.

<u>Proof of the corollary</u>: Let $\sigma$ be the signature with *k* plus signs and $k-1$ minus signs of the form $+\cdots+-\cdots-$. The minimum-length sequence is $\mu(\sigma) = 1,...,1, 0, (1,0), ..., (1,0), 1$ with $k$+1 leading 1's, followed by a 0, then $k$–2 pairs 1,0, and ending in a 1, which is therefore of length $\lambda(\mu) = (k+1)+1+2(k-2)+1 = 3k-1$. But as we will show in the next section, the length of the minimum-length sequence of a signature $\sigma$ does not depend on the permutation of + and – signs, only on $k(\sigma)$, so that for general signatures, $\lambda(\mu(\sigma)) = 3k(\sigma)-1$. Therefore *n*-sequences can only have signatures with $k \le \lfloor (n+1)/3 \rfloor$. For each such *k*, there are $\binom{2k-1}{k}$ permutations of + and – signs and for each of these, $m(n,\sigma) = n - \lambda(\mu(\sigma)) = n-(3k-1)$, which generate $\binom{m+k-1}{k-1} = \binom{n-(3k-1)+k-1}{k-1} = \binom{n-2k}{k-1}$ partitions by which to insert 0's before initial 1's. This yields the corollary. □

From Lemma 2 and Theorem 1, $h_4(n) = \sum_{i=2}^{n-1} h_2(i)$. This allows us easily to produce numerical tables such as the one below. In the original problem Litt posed with $n$=100 tosses, the excess number of wins for Bob over Alice is approximately $3.57382892 \times 10^{28}$, which is approximately 2.82% of the $2^{100}$ total number of tosses.

| $n$ | $h_2(n) = \Delta_n$ | $h_4(n) = D_n$ | $n$ | $h_2(n) = \Delta_n$ | $h_4(n) = D_n$ |
|---|---|---|---|---|---|
| 2 | 1 | 0 | 14 | 1,137 | 1,232 |
| 3 | 1 | 1 | 15 | 2,249 | 2,369 |
| 4 | 1 | 2 | 16 | 4,337 | 4,618 |
| 5 | 4 | 3 | 17 | 8,402 | 8,955 |
| 6 | 7 | 7 | 18 | 16,495 | 17,357 |
| 7 | 10 | 14 | 19 | 32,179 | 33,852 |
| 8 | 23 | 24 | 20 | 62,707 | 66,031 |
| 9 | 46 | 47 | 21 | 122,916 | 128,738 |
| 10 | 79 | 93 | 22 | 240,837 | 251,654 |
| 11 | 157 | 172 | 23 | 471,456 | 492,491 |
| 12 | 315 | 329 | 24 | 925,061 | 963,947 |
| 13 | 588 | 644 | 25 | 1,816,610 | 1,889,008 |

## 4. Formulas for the number of sequences with arbitrary score outcomes.

By counting arguments like those of Lemma 2, we arrive at the following formulas for *heady-s* binary $n$-sequences, i.e., sequences with $x_n = 1$ and $S(x) = s$, and *taily-s* binary $n$-sequences, i.e., sequences with $x_n = 0$ and $S(x) = s$, for any $s \in \mathbb{Z}$. Let $H_s(n)$ and $T_s(n)$ denote the number of heady-$s$ and taily-$s$ binary $n$-sequences, respectively.[1]

Theorem 2. For heady-$s$ binary $n$-sequences,

$$H_s(n) = \sum_{k=0\vee(-s)}^{\lfloor n_s/3 \rfloor} \binom{2k+s}{k}\binom{n_s - 2k}{k}, \quad \text{where } n_s = n - s - 1 = 0, 1, \ldots, \tag{4.1}$$

for any integer $s$ satisfying

$$-\lfloor (n-1)/2 \rfloor \leq s \leq n-1. \tag{4.2}$$

$H_s(n) = 0$ for $s$ outside of (4.2). For taily-$s$ binary $n$-sequences,

$$T_s(n) = 1_{s=0} + \sum_{k=1\vee(-s)}^{\lfloor n_s'/3 \rfloor} \binom{2k+s-1}{k-1}\binom{n_s' - 2k}{k}, \quad \text{where } n_s' = n - s = 0, 1, \ldots, \tag{4.3}$$

for any integer $s$ satisfying

$$-\lfloor n/2 \rfloor \leq s \leq 0 \vee (n-3). \tag{4.4}$$

$T_s(n) = 0$ for $s$ outside of (4.4).

[1] Note that what we had previously denoted by $h_2(n), h_3(n), h_4(n), t_2(n), t_3(n),$ and $t_4(n)$ would now be written $H_1(n), H_0(n), H_{-1}(n), T_1(n), T_0(n),$ and $T_{-1}(n)$. Cumulative frequencies like $h_1(n)$ would be sums such as $H_2(n) + H_3(n) + \cdots$.

Proof. For $H_s(n)$, first consider non-negative values of $s$, and for $k \geq 0$, let σ be any signature with $k$ minus signs and $k+s$ plus signs.[2] There are $\binom{2k+s}{k}$ rearrangements of σ with any given $k$. For the canonical signature $+\cdots+-\cdots-$, the minimum-length heady sequence is $m = 1,...,1,\ 1,0,...,1,0,\ 1$ with $k+s$ 1's followed by $k$ pairs 1,0, followed by a terminal 1, of total length $\lambda = (k+s)+2k+1 = 3k+s+1$. As we show below in Lemma 3, for any given signature σ, the length of the corresponding minimum-length sequence $m(\sigma)$ depends only on the number of plus and minus signs, so $\lambda(m(\sigma)) = 3k+s+1$ for σ or any of its rearrangements. Therefore, we must have $n \geq 3k+s+1$, i.e., $k \leq \left\lfloor \frac{n-s-1}{3} \right\rfloor = \lfloor n_s/3 \rfloor$. For a heady-$s$ binary $n$-sequence $x$ with the given signature, there are $n-\lambda(m(\sigma)) = n-s-1-3k = n_s - 3k$ additional zeros that may appear before the existing zeros in $m(\sigma)$, or before the initial 1 in $m(\sigma)$, i.e., in $k+1$ possible positions, and the number of multinomial partitions of $n-\lambda$ in $k+1$ cells is $\binom{n_s - 3k + (k+1) - 1}{(k+1)-1} = \binom{n_s - 2k}{k}$. Thus $H_s(n)$ equals (4.1) as claimed for non-negative $s$. However, nowhere in the argument did we use the assumption of non-negative $s$ apart from the tacit requirement that the number of plus signs $k+s$ should be non-negative. Thus (4.1) also applies for $s<0$, subject to the requirement that $k \geq -s$, in which case $\binom{2k+s}{k} > 0$ as well. Thus, the summation in (4.1) may start with $k = 0 \vee (-s)$ as in (4.1). For the range of $s$, the largest value of $s$ with positive $H_s(n)$ is plainly $s=n-1$, for the largest possible score results from $x=1,\dots,1$ with $S(x) = n-1$. Furthermore, for the summation in (4.1) to be non-vacuous, we must have $-s \leq n_s/3 = (n-s-1)/3$, i.e., $-2s \leq n-1$, i.e., $s \geq \lceil -(n-1)/2 \rceil = -\lfloor (n-1)/2 \rfloor$. This yields the range of $s$ specified in (4.2).

For $T_s(n)$, we note that any taily-$s$ binary $n$-sequence has a signature ending in a minus sign, except for $x=0,\dots,0$, which has a null signature and a zero score. Thus $T_s(n)$ equals $1_{s=0}$ plus the terms for non-null signatures. For those, let $k$ again denote the number of minus signs in signature σ with $k+s$ plus signs, this time with $k \geq 1$. There are now $\binom{2k+s-1}{k-1}$ rearrangements of σ since as noted, all must end with a minus sign. For the canonical signature $+\cdots+-\cdots-$, the minimum-length taily sequence is $m = 1,...,1,\ 1,0,...,1,0$ with $k+s$ 1's followed by $k$ pairs 1,0, of total length $\lambda = (k+s)+2k = 3k+s$, the same for any rearrangement σ. So we must have $n \geq 3k+s$, i.e., $k \leq \left\lfloor \frac{n-s}{3} \right\rfloor = \lfloor n_s'/3 \rfloor$. For any rearrangement σ there are $n-\lambda(m(\sigma)) = n-s-3k = n_s' - 3k$ additional zeros that may appear before the existing zeros in $m(\sigma)$, or before the initial 1 in $m(\sigma)$, i.e., again in $k+1$ possible positions, with $\binom{n'_s-3k+(k+1)-1}{(k+1)-1}=\binom{n'_s-2k}{k}$ partitions. Thus $T_s(n)$ equals $1_{s=0}$ plus a sum over $k$ of products $\binom{2k+s-1}{k-1}\binom{n'_s-2k}{k}$. Note that the largest value of $s$ for given $n$ is $s=n-3$, for the largest possible score results from taily $x=1,\ldots,1,0$ with score $S(x)=n-3$. So for non-vacuous summation, we must have $s\le n-3$, in which case $n'_s=n-s\ge 3$ and the upper limit of summation is $\lfloor n'_s/3\rfloor\ge 1$. For $n<s-3$, the summation is vacuous; in the special case $s=0$, $T_0(1)=1$ corresponding to the taily 1-sequence 0 with score 0, and $T_0(2)=1$ as well, corresponding to the taily 2-sequence 0,0. For negative $s$, non-zero terms in the summation begin with $k=-s\ge 1$, for which $\binom{2k+s-1}{k-1}>0$. Thus the summation may start with $k=1\vee(-s)$ for any $s$. Furthermore, the range of $k$ is vacuous unless $-s\le\lfloor n'_s/3\rfloor\le(n-s)/3$, i.e., $-2s\le n$, i.e., $s\ge\lceil -n/2\rceil=-\lfloor n/2\rfloor$. This yields (4.4). □

[2] In Lemma 2 we let $k$–1 denote the number of minus signs, but it will be more convenient to use $k$ for that purpose here.

**Remark.** By noting how heady-$s$ and taily-$s$ binary sequences arise, Theorem 2 generates a family of combinatorial identities that might not otherwise be obvious. For example, by noting that taily-$s$ $n$-sequences can only arise from taily $s$ $(n-1)$-sequences followed by a terminal 0 or from heady-$(s+1)$ $(n-1)$-sequences followed by a terminal 0, we get the identity $T_s(n)=T_s(n-1)+H_{s+1}(n-1)$. □

In the above arguments we used the invariance of the minimum length of heady-$s$ or taily-$s$ sequences under permutations of plus and minus signs of given signatures. We demonstrate that next.

Lemma 3. The length $\lambda(m(\sigma))$ of the minimum-length sequence $m(\sigma)$ for a given non-null signature $\sigma$ with $k$ minus signs and $k+s$ plus signs depends only on $k$ and $s$ and whether $m(\sigma)$ is heady-$s$ or taily-$s$.

Proof. Any sub-signature $\cdots-+\cdots$ of σ can be transposed to sub-signature $\cdots+-\cdots$ of σ′, say, with identical signs in both pre- and post-elided segments without changing the length of the minimum-length binary sequence. This is because in the first case, $m(\sigma)$ is of the form $\cdots 1,0,1,1,\cdots$ while in the second case, $m(\sigma')$ is of the form $\cdots 1,1,0,1,\cdots$, since the bit of $m(\sigma')$ following the displayed 0 must be a 1 whether the sign following the minus sign in $\sigma'$ is a + or a – sign. Thus $\lambda(m(\sigma))=\lambda(m(\sigma'))$. Since any signature can be rearranged to canonical form by a sequence of –+ to +– transpositions, $\lambda(m(\sigma))$ depends only on the number of plus and minus signs. □

## 5. Recursive algorithms for generating $H_s(n)$ and $T_s(n)$.

For generating tables of $H_s(n)$ and $T_s(n)$ for $n=1,\ldots,N_{max}$, it is more efficient to use recursive calculation rather than recalculating (4.1) or (4.3) for each $n$. We only know of one recursion of fixed order, for $H_1(n)$ (see Levin, 2024). In this section we state recursive algorithms of increasing order of intermediate stored arrays for arbitrary $s$. We use a hybrid APL notation for brevity.

Algorithm 1. Input: $N_{max}$. Output: nested arrays $H_s(n)$ for $n=1,\ldots,N_{max}$ and $-\lfloor (n-1)/2 \rfloor \le s \le n-1$.

Step (0). Initialize nested vector *Z* of length 2: *Z*[1]←⊂1 2ρ0 1 and *Z*[2]←⊂2 2ρ1 1 0 1. The elements of *Z*[*n*] give the pairs $s, t_s(n)$ for $n=1$ and 2. Initialize array *ts*←⊃*Z*[2] and 2×2 nested array *us*←2 2ρ1, (⊂1), 0, (⊂1). The first column of *us* contains $s=1$, 0 and the second column contains nested vectors $u_s(2)=1$ and 1, each of length 1, for $s=1$ and 0, respectively. $u_s(n)$ will contain $\binom{2k+s}{k}\binom{n_s-2k}{k}$ for relevant values of *k*. Also initialize array *vl*←2 2ρ1 1 0 1 containing $v_s(0)=1$ for $s=1$ and 0. $v_s(l)$ will contain $\binom{2l+s}{l}$ for relevant values of *l* equal to the maximum value of index *k*. Set *n*←3 and execute the following steps (1) through (6) until $n=N_{max}$.

(1) Loop on *s* from −⌊(*n*−1)÷2 to *n*−1. Go to step (6) when complete. Set *ns*←(¯2+*n*−*s*)÷3 and *l*←⌊*ns*÷3. Note that *l* and *ns* refer to the *previous* value of *n*, namely *n*−1.

(2) If *s* is not present in *ts*[;1], then append *ts*←*ts*,[1]*s*,1; *us*←*us*,[1]*s*,1; and *vs*←vs,[1]*s*,1. Go to next *s* in step (1). If *s* is present in *ts*[;1], then set *w* to the position of *s* in *ts*[;1] and set *t*←*ts*[*w*;2] and *u*←*us*[*w*;2].

(3) Set vector *k*←(0⌈−*s*),…,*l* and set *t*←*t*++/*inc*←*u*×*k*÷*ns*+1−3×*k* and *u*←*u*+*inc*, these being componentwise vector operations.

(4) If *ns* is not congruent to 2 modulo 3, go to step (5). If *ns* is congruent to 2 mod 3, then set *v*←*vs*[*w*;2] and *t*←*t*+*v*←*v*×(*s*+2+2×*l*)×(*s*+1×*l*)÷(*l*+1)×(*l*+*s*+1). Then append *u*←*u*,*v* and update *vs*[*w*;2] ←*v*.

(5) Update *ts*[*w*;2]←*t* and *us*[*w*;2]←⊂*u*, then go to next *s* in step (1).

(6) Update *Z*←*Z*,⊂*ts*←¯1ϕ[1]*ts*; *us*←¯1ϕ[1]*us*; and *vs*←¯1ϕ[1]*vs*. [These steps rotate the arrays circularly one row downward to arrange them in decreasing order of *s*.] Then go to the next *n*←*n*+1 and stop after completing the case $n=N_{max}$.

Algorithm 2. Input: $N_{max}$. Output: nested arrays $T_s(n)$ for $n=1,\ldots,N_{max}$ and $-\lfloor n/2 \rfloor \le s \le 0 \vee (n-3)$.

Step (0). Initialize nested vector *Z* of length 3: *Z*[1]←⊂1 2ρ0 1; *Z*[2]←⊂2 2ρ0 1 ¯1 1; and *Z*[3]←⊂2 2ρ0 2 ¯1 2. Each element of *Z*[*n*] gives the pair $s, t_s(n)$. Initialize array *ts*←⊃*Z*[3] and 2×2 nested array, *us*←2 2ρ0, (⊂1), ¯1, (⊂2). The first column of *us* contains $s=0,$¯1 and the second column contains nested vectors $u_s(3)=2$ and 2, each of length 1, for $s=0$ and ¯1, respectively. $u_s(n)$ will contain $\binom{2k+s-1}{k-1}\binom{n_s'-2k}{k}$ for relevant values of *k*. Also initialize array *vs*←2 2ρ0 1 ¯1 1 containing $v_s(3)=1$ for $s=0$ and ¯1. $v_s(l)$ will contain $\binom{2l+s-1}{l-1}$ for relevant values of *l* equal to the maximum value of index *k*. Set *n*←4 and execute the following steps (1) through (6) until $n=N_{max}$.

(1) Loop on *s* from −⌊*n*÷2 to *n*–3. Go to step (6) when complete. Set *ns*←(¯1+*n*–*s*)÷3 and *l*←⌊*n*÷3. Note that *l* and *ns* refer to the *previous* value of *n*, namely *n*–1.

(2) If *s* is not present in *ts*[;1], then append *ts*←*ts*,[1]*s*,1; *us*←*us*,[1]*s*,1; and *vs*←v*s*,[1]*s*,1. Go to next *s* in step (1). If *s* is present in *ts*[;1], then set *w* to the position of *s* in *ts*[;1] and set *t*←*ts*[*w*;2] and *u*←*us*[*w*;2].

(3) Set vector *k*←(1⌈–*s*),…,*l* and set *t*←*t*++/*inc*←*u*×*k*÷*ns*+1–3×*k* and *u*←*u*+*inc*, these being componentwise vector operations.

(4) If *n* is not congruent to 2 modulo 3, go to step (5). If *n* is congruent to 2 mod 3, then set *v*←*vs*[*w*;2] and *t*←*t*+*v*←*v*×(*s*+1+2×*l*)×(*s*+2×*l*)÷*l*×(*l*+*s*+1). Then append *u*←*u*,*v* and update *vs*[*w*;2] ←*v*.

(5) Update *ts*[*w*;2]←*t* and *us*[*w*;2]←⊂*u*, then go to next *s* in step (1).

(6) Update *Z*←*Z*,⊂*ts*←¯1ϕ[1]*ts*; *us*←¯1ϕ[1]*us*; and *vs*←¯1ϕ[1]*vs*. [These steps rotate the arrays circularly one row downward to arrange them in decreasing order of *s*.] Then go to the next *n*←*n*+1 and stop after completing the case $n=N_{max}$.

**References.**